\documentclass[12pt,a4paper]{article}
\usepackage{graphicx,amsmath,amsfonts,amssymb,color,mathrsfs, amsmath,amsthm}
\usepackage{epsfig}

\usepackage{dsfont}

\newcommand{\chuhao}{\fontsize{19pt}{\baselineskip}\selectfont}

\newcommand{\BOX}{\hfill $\Box$}

\newcommand{\stsp}{\mathbb{E}}
\renewcommand{\P}{{\rm P}}
\newcommand{\E}{{\rm E}}
\newcommand{\ud}{\mathrm{d}}
\newcommand{\DD}{\mathcal{D}}

\newcommand{\indic}{1{\hskip -2.5 pt}\hbox{I}}

 \newcommand{\vc}[1]{{\boldsymbol #1}}
    \newcommand{\valpha}{\vc\alpha}
    \newcommand{\vmu}{\vc\mu}
    \newcommand{\vpi}{\vc\pi}
    \newcommand{\vtau}{\vc\tau}
    \newcommand{\vone}{\vc 1}
    \newcommand{\vzero}{\vc 0}

 \newtheorem{theorem}{Theorem}[section]
 \newtheorem{lemma}[theorem]{Lemma}
 
 \newtheorem{proposition}[theorem]{Proposition}
 
 \newtheorem{remark}[theorem]{Remark}

 \newtheorem{corollary}[theorem]{Corollary}

\newcommand{\vligne}[1]{\begin{bmatrix} #1 \end{bmatrix}}

\usepackage{pifont}

\title{\bf\color{black} \chuhao{ Poisson's equation for discrete-time quasi-birth-and-death processes}}
\author{
Sarah Dendievel\thanks{%
Universit\'e Libre de Bruxelles, D\'epartement d'Informatique,
CP~212, Boulevard du Triomphe, 1050 Bruxelles, Belgium; Sarah.Dendievel@ulb.ac.be, latouche@ulb.ac.be.}
\and Guy Latouche$^*$
\and Yuanyuan Liu\thanks{%
Corresponding author. School of Mathematics, Railway
Campus, Central South University, Changsha, Hunan, 410075, China; liuyy@csu.edu.cn.}
}

\date{This is an Author¡¯s Accepted Manuscript of an article published in Performance Evaluation, available online at:
http://dx.doi.org/10.1016/j.peva.2013.05.008}

\begin{document}
 \maketitle

\begin{abstract}
We consider Poisson's equation for quasi-birth-and-death processes (QBDs) and
we exploit the special transition structure of QBDs to obtain its solutions in two different forms.
  One is based on a decomposition through first passage times to lower levels, the other is based
  on a recursive expression for the deviation matrix.

We revisit the link between a solution of Poisson's equation  and
perturbation analysis and we show that it applies to QBDs.  We
conclude with the PH/M/1 queue as an illustrative example, and we
measure the sensitivity of the expected queue size to the initial value.
\medskip

\end{abstract}

\section{Introduction}

  Poisson's equation has the following form:
\begin{equation}
  \label{Poisson's equation}
  (I-P) {\vc h} = {\vc{g}},
\end{equation}
where $P$ is the transition matrix of a Markov chain on some denumerable state space  $\mathbb{E}$ and
     ${\boldsymbol g}$  is a given  vector on $\mathbb{E}$, subject to some constraints.
Equations of the form (\ref{Poisson's equation}) frequently occur in the analysis of Markov chains.  In particular, as  remarked  in Meyn and Tweedie~\cite[Pages 458-459]{mt09}, we find them in the context of  central limit theorems, perturbation theory, controlled Markov processes, variance
analysis of simulation algorithms, etc.  It is of interest to note that the equation has been analyzed in cases where the state space is continuous: Glynn \cite{glynn94} derives a solution of Poisson's equation for the waiting time process of the recurrent discrete-time M/G/1 queue, and Asmussen and Bladt \cite{ab94} consider the waiting time for continuous-time queues driven by a Markovian marked point process.

Here, we assume that the Markov chain is a quasi-birth-and-death
process (QBD): the state space is $\mathbb{E}=\{\cup_{i=0}^\infty
\ell(i)\}$, where $\ell(i)=\{(i,j): 1\leq j\leq m\}$ denotes the
level set, and the transition matrix is
\begin{equation}
   \label{e:P}
 P=
\begin{bmatrix}
   B & A_1 & 0&  \ldots \\
   A_{-1} & A_0 & A_1 &\ddots \\
   0 & A_{-1} & A_0 & \ddots \\
   \vdots & \ddots& \ddots & \ddots
\end{bmatrix}
   \end{equation}
where $A_{-1}$, $A_0$, $A_1$, and $B$ are square matrices of order
$m$. We assume  that the process is irreducible and that
the matrix $A$ defined as $A=A_{-1}+A_0+A_1$ is irreducible.  We
also assume that $m < \infty$.  For details, we refer to
Neuts~\cite{neuts81} and Latouche and Ramaswami~\cite{lr99}.

We focus on the case where the QBD is positive recurrent with  invariant distribution
${\vc \pi}$.  Such is the case if and only if
 $  \boldsymbol{\mu} (A_{-1} -A_1) \vc{1} > 0$,
 where $\vmu$ is the invariant probability vector of $A$, and  ${\boldsymbol 1}$ is a
 column vector composed of~$1$s.  Equivalent drift conditions with different forms may be found in
  Latouche and Taylor~\cite{{lt00}}.
 The  distribution $\vc{\pi}$ is decomposed as $\vc\pi =
 \vligne{\vpi_0 & \vpi_1 & \vpi_2 & \ldots}$ where $\vpi_n$ is the
 sub-vector of stationary probability for level $n$ and we similarly
 decompose the vector $\vc g$ into sub-vectors.

Our main objective is to use structural properties of QBDs to
express the solutions of Poisson's equation. In
Section~\ref{s:basicPoisson} we give some basic properties about the
 equation for general countable Markov chains, we investigate
the relation between the solution $\vc h$ and the deviation matrix of $P$, and we
obtain a decomposition of the system~(\ref{Poisson's equation}) through
a first return time argument.  This is the key to our result in
Section~\ref{s:returnTimes} where $P$ has the structure~(\ref{e:P}), and we
obtain a solution of (\ref{Poisson's equation}) in terms of first passage times to lower levels.
In Section~\ref{s:deviationMatrix}, an explicit expression of the
deviation matrix for QBDs is obtained, which is of interests in
its own right.

It is shown in Liu and Hou~\cite{Liu2006} that a positive recurrent
QBD is geometrically ergodic and we determine an explicit drift
condition in Section~\ref{s:returnTimes}.  Motivated by that
observation, we revisit a few well-known properties of the solutions
of Poisson's equation, proved earlier for general Markov chains under
more restrictive assumptions,
and we show that they hold in the case of QBDs.
The link between the solution and sensitivity analysis is investigated
in Section \ref{s:perturbation} --- our results improve on the
corresponding ones of Cao and Chen~\cite{cc97}.  Finally, we give an
illustrative example in Section~\ref{s:exemple}: we compute the total
difference between the expected number of customers of a PH/M/1 queue
in finite time and its stationary value.

\section{General Markov chains}
   \label{s:basicPoisson}

 Let $\{X_t\}$ be a
discrete-time Markov chain, irreducible, aperiodic and positive recurrent  on a countable state space
$\mathbb{E}$, with transition matrix $P$ and
invariant distribution ${\vc \pi}$.
Poisson's equation is defined in Makowski and Schwartz~\cite{ms02} as
\begin{equation}\label{2.1}
   ( I-P) \boldsymbol{h}=\boldsymbol{g}-\omega\boldsymbol{1},
\end{equation}
where $ \boldsymbol{g}$ is a given column vector, the scalar $\omega$
and the vector $\vc h$ forming together the solution of (\ref{2.1}).

If the state space $\stsp$ is finite, then the solutions are
\begin{equation}
   \label{e:drazin}
\vc h = (I-P)^{\#} {\vc{g}} + c \vone   \qquad \mbox{and} \qquad \omega = \vpi \vc g,
\end{equation}
where  $(I-P)^{\#}$ is the group inverse of $I-P$ (see Meyer~\cite{meyer75}, Campbell and Meyer~\cite{cm91}),
and $c$ is an arbitrary constant.
Actually, if $\stsp$ is finite, it is obvious  that $\omega$ has to be
equal to $\vpi \vc
g$: with $\vc\pi (I-P) = \vzero$, there cannot be any other solution of (\ref{2.1}).  Thus,
we might have written the system  as $(I-P) \vc h = \overline{\vc{g}}$, where
$\overline{\vc{g}} = \vc g - (\vpi \vc g) \vone$.

One simple interpretation of (\ref{e:drazin}) goes as follows.  Assume
that $\vc g$ is a vector of state-dependent rewards: the reward at
time $t$ is $g_i$ if $X_t=i$.  The stationary expected reward per unit
of time is $\vc\pi \vc g$ and if $\vpi\vc g = 0$, one shows that the
$i$th component of $(I-P)^\# \vc g$ is the total expected reward
accumulated over the whole history, given that $X_0=i$.   If $\vpi\vc
g \not= 0$, then the total reward diverges to $+\infty$ or $-\infty$, and
$(I-P)^\# \vc g$ is the vector of total expected {\em difference}
between the actual reward and its stationary mean, given the initial
state.  The arbitrary constant $c$ in (\ref{e:drazin}) reflects the
fact that $I-P$ is singular, and that an additional constraint is
needed to thoroughly specify the solution of Poisson's equation.

If $\stsp$  is infinite, the situation is more involved but one does
have the following property.
It is a direct consequence of   \cite[Theorem 9.5]{ms02} and we only
need to verify that the assumptions there  both hold.  Defining $T(j)$
to be the first return time to some state $j$,
the recurrence condition requires that $\E[ T(j)|X_0=i]$ be finite for
all $i$.  This results from the assumption that the Markov chain is
irreducible and positive recurrent.  The integrability condition
requires that  $\E[\sum_{0 \leq t < T(j)} |g_{X_t}|\,|X_0=i]$  be
finite for all $i$.  This is a consequence of~\cite[Theorem~14.1.3]{mt09}.

\begin{lemma}
   \label{t:basicsolution}
Assume that the Markov chain $\{X_t\}$ is irreducible and positive
recurrent, and that $\vpi |\vc g| < \infty$. Take $j$ to be a fixed,
arbitrary state and define the first {return time} to $j$ as $T(j) =
\inf\{t \geq 1: X_t =j\}$; also, define
\[
\zeta_i = \E[\sum_{0 \leq t < T(j)} g_{X_t}|X_0=i]
  \qquad \mbox{and} \qquad
 \tau_i = \E[ T(j)|X_0=i]
\]
for $i \in \stsp$.
The pair $(\vc h, \omega)$ given by
\begin{equation}
   \label{e:solution}
\omega = \vpi \vc g
\qquad \mbox{and} \qquad
h_i = \zeta_i- \omega \tau_i, \qquad i \in \stsp
\end{equation}
is one solution of (\ref{2.1}) with $h_j=0$.
\BOX
\end{lemma}

\begin{remark} \rm
In the sections to follow, we restrict ourselves to solutions which are based on first return times,
as exposed in the lemma.  If $\vpi |\vc g| < \infty$, we may rewrite (\ref{e:solution}) as
\[
h_i = \E[\sum_{0 \leq t < T(j)} (g_{X_t} - \omega) |X_0=i] =
\E[\sum_{0 \leq t < T(j)} \overline{{g}}_{X_t} |X_0=i]
\]
where $\overline{\vc{g}} = \vc g - \omega \vone$ is such that $\vpi \overline{\vc{g}} =0$.
 This is the reason why we assume later that $\vpi \vc g =0$.
\end{remark}

\begin{remark} \rm
In the statement of Lemma~\ref{t:basicsolution}, the state $j$ is arbitrary.
Thus, the lemma defines a different solution for each choice of $j$.
That choice, however, is mostly irrelevant as  one shows that any two such solutions differ by a constant.
In general, uniqueness of the solution of (\ref{2.1}), up to an additive constant,
 is not readily obtained when $\stsp$ is infinite, unless there
 is some integrability constraint imposed on the solution, see \cite[Section 9.7]{ms02} for a detailed example.
\end{remark}

We now extend Lemma~\ref{t:basicsolution} and express the solution of Poisson's equation
 in terms of first return times to {\em subsets} of states.  As we show  in the proof of Theorem~\ref{th2.1},
  there are circumstances where it is more convenient to work with return times to a set of states than to a single state.


\begin{lemma}
   \label{t:returntoA}
Assume that the Markov chain $\{X_t\}$ is irreducible and positive
recurrent, and that $\vpi |\vc g| < \infty$.
 Take $A$ to be an arbitrary, non-empty, subset of $\stsp$,  define  $T(A) = \inf\{t \geq 1: X_t \in A\}$
  to be the first return time to $A$, and define
for all $i \in \stsp$
\[
y_i(A) = \E[\sum_{0 \leq t < T(A)} g_{X_t}|X_0=i]
  \qquad \mbox{and} \qquad
 \tau_i(A) = \E[ T(A)|X_0=i].
\]
The solution $\boldsymbol{h}$ defined in Lemma~\ref{t:basicsolution}
satisfies
\begin{equation}\label{2.1.1}
h_{i} = y_i(A) - \omega \tau_i(A) +
\E[ h_{X_{T(A)}}| X_0 = i], \qquad i\in \mathbb{E}.
\end{equation}
\end{lemma}
\proof
Take an arbitrary state $j$ in $A$ and take the vector $\vc h$ defined in (\ref{e:solution}).
 Since $T(j)=\max(T(j), T(A))$, we may write
\begin{align*}
h_i &= \E[\sum_{0 \leq t < \max(T(j),T(A))} g_{X_t}|X_0=i]  - \omega \E[\max( T(j),T(A))|X_0=i] \\
 & = \E[\sum_{0 \leq t < T(A)} g_{X_t}|X_0=i]   + \E[\sum_{T(A) \leq t < T(j)} g_{X_t}|X_0=i] \\
& \quad - \omega \E[T(A)|X_0=i]     - \omega \E[( T(j)-T(A)) \indic{[T(j)>T(A)]} |X_0=i],
\end{align*}
where $\indic{[\cdot]}$ is the indicator function. By the strong Markov property, this proves the claim.
\BOX

Let us  partition the state space as $\mathbb{E}=A\cup B$, where $A$ and $B$ are two proper subsets,
and partition in a similar manner the transition matrix $P$ as
\begin{equation}\label{par}
P=
\begin{bmatrix}
P_{A} & P_{AB} \\
P_{BA} & P_{B}
\end{bmatrix}
\end{equation}
and the vector $\boldsymbol{h}$ as
\[
\vc h = \vligne{\boldsymbol{h}_A \\
  \boldsymbol{h}_B}.
\]
 The next theorem will be useful in solving Poisson's equation for QBDs.

\begin{theorem}\label{th2.1}
Assume that the Markov chain $\{X_t\}$ is irreducible and positive
recurrent, and that $\vpi |\vc g| < \infty$. The vector
$\boldsymbol{h}$ given by
\begin{equation}\label{2.1.5}
\boldsymbol{h}_{A}=\boldsymbol{y}_A(A)
-( \boldsymbol{\pi g}) \boldsymbol{\tau}_{A}(A)
+(
P_{A}+P_{AB}N_BP_{BA}) \boldsymbol{h}_{A}
\end{equation}
and
\begin{equation}\label{2.1.6}
\boldsymbol{h}_{B}=\boldsymbol{y}_{B}(A)
-( \boldsymbol{\pi g}) \boldsymbol{\tau}_{B}(A) +  N_B P_{BA} \boldsymbol{h}_{A}
\end{equation}
is a solution of the Poisson equation, where $N_B=\sum_{n=0}^\infty P_B^n$.
\end{theorem}

\proof
We rewrite (\ref{2.1.1}) as
\[
h_i = y_i(A) - \omega \tau_i(A) + \sum_{k \in A} \P[X_{T(A)} = k|X_0=i] h_k(A).
\]
Since
\begin{align*}
\P[X_{T(A)} = k|X_0=i] & = (\sum_{n=0}^\infty P_B^n P_{BA})_{ik} && \mbox{for $i \not\in A$,}  \\
  & = (P_{A}+P_{AB} \sum_{n=0}^\infty P_B^n P_{BA})_{ik}  && \mbox{for $i \in A$,}
\end{align*}
we observe that the system  (\ref{2.1.1})  is equivalent to (\ref{2.1.5},
\ref{2.1.6}), and
the theorem follows.
\BOX

To conclude this section, we obtain an expression similar to (\ref{e:drazin}) in the case where $\stsp$ is infinite.
The natural extension of the group inverse to infinite-sized aperiodic stochastic matrices is the {\em deviation matrix}
 $\DD$ defined as
\begin{equation}
   \label{e:deviation}
\DD = \sum_{n \geq 0} (P^n - \vone \vpi),
\end{equation}
assuming that it exists.  If $\stsp$ is finite and $P$ is irreducible, then $\DD=(I-P)^\#$
by Campbell and Meyer~\cite[Theorem 8.3.1]{cm91}.

Denoting by $\E_\pi[\cdot]$ the conditional expectation given that $X_0$ has the distribution $\vpi$,
 we have the following criteria for the existence of $\DD$.

\begin{lemma}
   \label{t:existence}
Assume that the Markov chain $\{X_t\}$ is irreducible, aperiodic and
positive recurrent. The deviation matrix $\DD$ exists if and only if
$\E_\pi[T(j)] < \infty$ for some $j$ in $\stsp$ or, equivalently, if
$\E[T^2(j)|X_0=j] < \infty$.  If the property holds for one state,
then $\E_\pi[T(i)] < \infty$ for all $i$ in $\stsp$.
\end{lemma}

\proof
This is nearly a restatement of Syski~\cite[Proposition~3.2]{syski78}, the difference being that the analysis
 in \cite{syski78} is for continuous-time Markov chains.

Starting from the discrete-time Markov chain $\{X_n\}$, we may define a continuous-time Markov chain $\{X(t)\}$
 by making  the process change state at the epochs of transition of a Poisson process with constant parameter.
   As shown in Coolen-Schrijver and van Doorn~\cite[Section 2]{cv02}, the deviation matrix of the discrete-time
    chain exists if and only if the deviation matrix of the continuous-time chain exists.
\BOX

When $\DD$ exists, its elements are given by
\begin{eqnarray}
   \label{e:djj}
\DD_{jj}  & = & \pi_j (\E_\pi[T(j)] -1) \\
   \label{e:dij}
\DD_{ij}  & = & \DD_{jj} - \pi_j \E[T(j)|X_0=i], \qquad  \mbox{for
  $i$, $j$ in $\stsp$, $i \not= j$}
\end{eqnarray}
(see \cite{cv02, syski78}).  The connection with Poisson's equation is established by the following two lemmas.

\begin{lemma}
   \label{t:DandPa}
If the Markov chain $\{X_t\}$ is irreducible, aperiodic and positive
recurrent, and if $\E_\pi[T(j)] < \infty$
 for some $j$ in $\stsp$, then $\DD$ is the unique solution of the system
\[
(I-P)\DD = I-\vone \vpi, \qquad  \vpi \DD = \vzero
\]
in the set of matrices $M$ such that $\vpi |M|$ is finite.
Furthermore, $\DD\vone = \vzero$.
\end{lemma}
\proof
The fact that $\DD$ has the stated properties is proved in \cite[Theorem 5.2]{cv02}
 so that we only need to prove uniqueness.  Recall that the finiteness of $\E_\pi[T(j)]$
 for one state $j$ implies the finiteness of $\E_\pi[T(k)]$ for all state $k$.

Take an arbitrary state $k$ and define $\vc d$ as the $k$th column of $\DD$.
By (\ref{e:djj}, \ref{e:dij}), we may write that $\vc d = \DD_{kk} \vone - \pi_k \vc\tau^*(k)$,
 where $\tau^*_i(k) = \E[T(k)|X_0=i]$ for $i \not= k$ and  $\tau^*_k(k) = 0$.   Now,
\[
\E_\pi[|d_{X_n}|] = \sum_{i \in \stsp} \pi_i |d_i| \leq |\DD_{kk}|
+ \vpi |\vc\tau^*(k)| \leq |\DD_{kk}| + \E_\pi[T(k)] < \infty
\]
for all $n$.  By \cite[Theorem 9.1]{ms02},  all matrices $X$ solution
of $(I-P)X= I - \vone \vpi$ and such that $\vpi |X|$ is finite are given by
 $X = \DD + \vone \vc x$ for some vector $\vc x$.  Since, in addition,
 $\vpi X$ must be equal to $\vzero$, we conclude that $\vc x = \vzero$ and that $X = \DD$.
\BOX

Finally, we define the $w$-norm $\| \cdot \|_w$ as
\[
\| A \|_w  = \sup_{i \in \stsp} (1/w_i) \sum_{j \in \stsp} |A_{ij}| w_j,
\]
where the weight function $\vc w$ is bounded away from zero:
$w_i < \infty$ for all $i$, $\inf_{i \in \stsp} w_i >0 $.
A Markov chain is $w$-geometrically ergodic if $\|P\|_w < \infty$ and
if there are scalars $r >0$ and $\beta < 1$ such that
$\|P^n - \vone \vpi \|_w \leq r \beta^n$ for $n \geq 0$.   With this,
Bhulai and Spieksma~\cite[Theorem 4]{bs03} implies the following theorem.

\begin{theorem}
   \label{t:DandPb}
Assume that the chain $\{X_t\}$ is irreducible, aperiodic and positive
recurrent, that $\vpi |\vc g|< \infty$ and that $\vpi \vc g =0$.

If the Markov chain $\{X_t\}$ is $w$-geometrically ergodic for some $\vc
w$, then a solution of the Poisson equation (\ref{Poisson's
equation}) is given by $\vc h = \DD \vc g +c \vone$, where $c$ is
arbitrary. \BOX
\end{theorem}

Furthermore, one defines the drift condition \textbf{D}($\vc v, \lambda, b, C$)
 as follows: there exists a finite drift function $\boldsymbol{v}$ bounded away from zero,
  two constants $\lambda\in (0, 1)$ and $b<\infty$, and some finite set of states $C\subset \mathbb{E}$ such that
\begin{equation*}
  P\boldsymbol{v} \leq \lambda \boldsymbol{v} + b \indic{[{C}]}.
\end{equation*}
It follows from \cite[Theorem 16.01]{mt09} that a Markov chain is $w$-geometrically ergodic
if and only if the drift condition \textbf{D}($\vc v, \lambda, b, C$) holds for a vector $\boldsymbol{v}$,
 equivalent to $\boldsymbol{w}$ in the sense that for some $c>1$, $c^{-1} \vc w\leq \boldsymbol{v}\leq c\boldsymbol{w}$.

\section{The case of QBDs}
   \label{s:returnTimes}

Liu and Hou~\cite[Remark 3.2]{Liu2006} show that a positive
recurrent, aperiodic,  QBD is geometrically ergodic, so that
{Theorem~\ref{t:DandPb}} applies.  It is of independent interest to
determine a suitable set of parameters for \textbf{D}($\vc v,
\lambda, b, C$), and this we do now.

We define the matrix function $A(z)=\frac{1}{z}A_{-1}+A_0  + A_1 z$,
for $z > 0$, and we denote by $\sigma(z)$ the Perron-Frobenius
eigenvalue of $A(z)$ and by $\vc u(z)$ the corresponding strictly
positive right-eigenvector.   The eigenvector $\vc u(z)$ and the
scalar $\sigma(z)$ play a role in the theorem below.

\begin{lemma}
   \label{t:drift}
If the QBD with transition matrix (\ref{e:P}) is irreducible and positive
recurrent, then
the drift condition
\textbf{D}($\boldsymbol{v}, \lambda_0, b, C$) holds, with $\vc v_i=z_0^i
\boldsymbol{u}(z_0)$, $\lambda_0= \sigma(z_0)$, $C=\ell(0)$, and
\begin{equation}
   \label{e:b}
b= \max_{1\leq j \leq m} (B \boldsymbol{v}_0 + A_1 \boldsymbol{v}_1- \lambda_0 \boldsymbol{v}_0)_j,
\end{equation}
where $\vc u(z)$, $\sigma(z)$  are defined above and $z_0$ is the
minimal solution of the equation $\sigma'(z)=0$ with $z>1$.
\end{lemma}
\proof Obviously, $\sigma(1)=1$.   Furthermore, $\sigma(z)$
is analytic on $z>0$.  Finally, one verifies that
$\sigma'(1)=\boldsymbol{\mu} (A_{1}-A_{-1}) \boldsymbol{1}$ and this
is strictly negative for positive recurrent QBDs.   Altogether, this
proves that  there exists $z>1$ such that
$\sigma(z)<1$.  Take $\boldsymbol{v}_i= z^i \boldsymbol{u}(z)$ and
$\lambda= \sigma(z)$ for any such $z$. For $i\geq 1$,
\begin{eqnarray*}
  (P\boldsymbol{v})_i  & =& A_{-1} \boldsymbol{v}_{i-1} +A_0 \boldsymbol{v}_{i} + A_1 \boldsymbol{v}_{i+1}  \\
  & =&  z^{i}( \frac{1}{z} A_{-1}+  A_0  +  z  A_1) \boldsymbol{u}(z)
  \  = \ \lambda  \boldsymbol{v}_i
\end{eqnarray*}
and for $i=0$,
\[
(P \boldsymbol{v})_0= B \boldsymbol{v}_0 + A_1 \boldsymbol{v}_1 \leq
  \lambda \boldsymbol{v}_0 +  (\max_{1\leq j \leq m} (B
  \boldsymbol{v}_0 + A_1 \boldsymbol{v}_1- \lambda
  \boldsymbol{v}_0)_j) \textbf{1}.
\]
The matrix
polynomial analyzed in Bean {\em et al.}~\cite{bblppt97} is
identical to $z A(1/z)$ and we know from \cite[Theorem 5]{bblppt97}
that there exists some $z_0$, $1< z_0 < \infty$, such that
$\sigma(z)$ is minimal for $z=z_0$.   Actually, the results in
\cite{bblppt97}  indicate that either  this $z_0$ is unique, or
there exists some interval over which $\sigma(z)$ is minimal and
constant, in the latter case, we choose  $z_0$ to be the minimal
solution of the equation
$\sigma'(z) =0$ with $z>1$..
\BOX

It is clear from the proof that the choice of $z$ is to some extent
arbitrary.   By taking $z=z_0$, we minimize $\lambda$ within our
construction.

\begin{remark} \em
Mao {\em et al.}~\cite{mtzz12} analyze the more general Markov chains of GI/G/1 type.
 Our characterization in Lemma~\ref{t:drift} for QBD processes is more explicit than the
  one in~\cite[Theorem 3.1]{mtzz12}.
\end{remark}

Several properties of QBDs are related to three key matrices, $R$, $G$, and $U$,
which are characterized as follows:
\begin{align}
   \label{e:r}
R_{ij} & = \E[\sum_{0 \leq t < T(\ell(n))} \indic{[X_t=(n+1,j)]}|X_0=(n,i)] \\
   \label{e:g}
G_{ij} & = \P[T(\ell(n))< \infty, X_{T(\ell(n))}=(n,j)| X_0= (n+1,i)],  \\
\intertext{for $n \geq 0$, and}
   \label{e:u}
U_{ij} & = \P[T(\ell(n))< T(\ell(n-1)), X_{T(\ell(n))}=(n,j)| X_0= (n,i)],
\end{align}
for $n \geq 1$.
We refer to Latouche and Ramaswami~\cite[Chapter 6]{lr99} for details.
Here, we only mention that all three matrices may be efficiently
computed, and that if the QBD is positive recurrent, then $G$ is
stochastic, $U$ is sub-stochastic, and the spectral radius of $R$ is
strictly less than 1.  Furthermore, the stationary distribution is
given by $\vpi_n = \vpi_0 R^n$, for $n \geq 0$, and $\vpi_0$ is the
unique solution of the system $\vpi_0 P_*= \vpi_0$, with $\vpi_0
(I-R)^{-1} \vone = 1$, where $P_*=B+A_1 G$.  Finally, a useful
relation is that $R=A_1 (I-U)^{-1}$.

Not surprisingly, the three matrices also play a key role in determining a solution
of Poisson's equation as we show next.  We first apply Theorem~\ref{th2.1} and,
to that end, we choose $A$ to be the level 0 
and $B$ to be the set of all other levels.  To simplify the notations,
we write $T$, $\vc y$ and $\vc\tau$, respectively, for $T(\ell(0))$, $\vc y(\ell(0))$ and $\vc\tau(\ell(0))$.

\begin{theorem}
   \label{t:qbda}
If the QBD is irreducible and positive recurrent, and if the vector
$\vc g$ satisfies
\begin{equation}
   \label{e:vectorg}
\sum_{n \geq 0} \vpi_0 R^n |\vc g_n| < \infty,
\qquad
\sum_{n \geq 0} \vpi_0 R^n \vc g_n = 0,
\end{equation}
then a solution of the Poisson equation $(I-P) \vc h = \vc g$  is given by
\begin{eqnarray}
 \boldsymbol{h}_{0} &=& (I - P_*)^{\#} \boldsymbol{y}_{0} + c \boldsymbol{1}, \label{3.1} \\
 \boldsymbol{h}_{n} &=& \boldsymbol{y}_{n}
  + G^{n}\boldsymbol{h}_{0}, \ \ n\geq 1, \label{3.2}
\end{eqnarray}
where  $\boldsymbol{y}$ is defined in Lemma 2.4, $P_*  = B+A_{1}G$,
and $c$  is an arbitrary constant.  Moreover,  the vector $\vc y$ is
explicitly given by
\begin{align}
    \boldsymbol{y}_{n} &= \sum_{0 \leq i \leq n-1} G^{i}(I-U)^{-1}
    \sum_{l \geq 0}  R^l \vc g_{n-i+l},  \qquad \mbox{for  $n\geq 1$,} \label{3.6} \\
     \boldsymbol{y}_0 &= \vc g_{0} + A_1 \boldsymbol{y}_1. \label{3.5}
\end{align}
\end{theorem}
\proof The physical meaning of the matrix $N_B P_{BA}$ in Theorem
\ref{th2.1} is that each entry $((n,i),(0,j))$
 is the probability that, starting from $(n,i)$, the QBD reaches level 0 in finite time, and
 that the first state visited there is $(0,j)$.  For $n=1$, this probability is $G_{ij}$,
 by (\ref{e:g}).  For $n >1$, the process must successively visit the levels $n-1$,
 $n-2$, \ldots, $1$ because of the skip-free structure of (\ref{e:P}), and each step down
 is controlled by the same transition matrix $G$.  Thus, $(N_B P_{BA})_{(n,i),(0,j)}=(G^n)_{ij}$, and
\begin{equation}
   \label{e:wam1}
N_B P_{BA}=
\begin{bmatrix}
G \\
G^2 \\
G^3 \\
\vdots
\end{bmatrix}.
\end{equation}
As we assume that $\vpi \vc g = 0$,  (\ref{3.2}) directly result from (\ref{2.1.6}).

On the other hand,
$P_A+P_{AB}N_B P_{BA}= B + A_1 G = P_*$ and
 (\ref{2.1.5}) becomes
$
\boldsymbol{h}_{0} = \boldsymbol{y}_{0}
 + P_* \boldsymbol{h}_{0}.
$
By (\ref{e:drazin}),
\[
 \boldsymbol{h}_{0} = (I - P_*)^{\#} \boldsymbol{y}_{0}  + c \boldsymbol{1}
\]
for some $c$, since $P_*$ is stochastic and $m< \infty$.

As observed above, the first passage time to level 0 is a sum of
first passage times from one level to the one immediately below, and
we may write $T=\theta_n + \theta_{n-1} + \cdots + \theta_1$ if the
initial state is in level $n$, where $\theta_k$ is the first passage
time from level $k$ to level $k-1$. If we interpret the vector $\vc
g$ as a vector of rewards associated with visits in the different
states, as we suggested at the beginning of Section~\ref{s:basicPoisson},
then $\vc y_n$ is the expected reward accumulated over the
trajectory from level $n$ to level $0$ and we may write
\begin{equation}
   \label{e:yn}
\vc y_n = \vc u_n + G \vc u_{n-1} + G^2 \vc u_{n-2} + \cdots + G^{n-1} \vc u_1,
\end{equation}
where $\vc u_k$ is the vector of accumulated rewards during the first passage time from level $k$ to level $k-1$:
\[
(\vc u_k)_i = \E [ \sum_{0 \leq t < \theta_k} g_{X_t}|X_0 = (k,i)].
\]
Now, during the first passage time $\theta_n$, the process may visit any number of times the states in higher levels and, by \cite[Remark 6.2.8]{lr99}, we may interpret $(R^l)_{ij}$ as the expected number of visits to $(n+l,j)$ before the first return to level $n$ given that the process starts at $(n,i)$, for all $l \geq 0$.  Thus, if we decompose the trajectory by the first return to level $n$, we find that
\[
\vc u_n = \sum_{l \geq 0} R^l \vc g_{n+l}  + U \vc u_n = (I-U)^{-1} \sum_{l \geq 0} R^l \vc g_{n+l}.
\]
This, together with (\ref{e:yn}), proves (\ref{3.6}).
To prove
(\ref{3.5}), we condition on the first transition of the Markov chain.
\BOX

\begin{corollary}
   \label{t:tau}
If the QBD is irreducible and positive recurrent, then the expected
first passage times from level $n$ to level 0 are given by
\begin{align}
   \label{3.3}
    \boldsymbol{\tau}_{0}&=  (I-R)^{-1} \boldsymbol{1} \\
\label{3.4}
     \boldsymbol{\tau}_{n}&= ((I-G^n)(I-G)^{\#} + n\boldsymbol{1 \gamma})(I-U)^{-1} (I-R)^{-1}\boldsymbol{1},
\end{align}
for $n \geq 1$,
where $\boldsymbol{\gamma}$ is the invariant probability vector of the stochastic matrix~$G$.
\end{corollary}
\proof
To prove this, we replace in (\ref{3.6}) the vector $\vc g$ by $\vone$.
This gives us $\vtau_1 = (I-U)^{-1}(I-R)^{-1} \vone$ and in general $
\vtau_n = (I + G + \cdots + G^{n-1}) \vtau_1$.  Since $\sum_{0 \leq i \leq n-1} G^i
= (I-G^n)(I-G)^{\#} + n \vone \vc\gamma$, this proves (\ref{3.4}).

Using the fact that $A_1 = R(I-U)$ (see [14, Page 137]), we obtain
from (\ref{3.5}) that
\begin{align*}
\vtau_0 & =\boldsymbol{1} + A_1\boldsymbol{\tau_1} = \vone + A_1 (I-U)^{-1} (I-R)^{-1} \vone \\
& = \vone + R (I-R)^{-1} \vone = (I-R)^{-1} \vone,
\end{align*}
which completes the proof.
\BOX

\section{Deviation matrix }
   \label{s:deviationMatrix}

It results from Lemma~\ref{t:DandPa} that the deviation matrix is a solution of
Poisson's equation and might be computed by applying Theorem~\ref{t:qbda}.
We give here another expression which is of independent interest.
To that end, we need the following preliminary result.  Define
\[
H = \begin{bmatrix}
A_0 & A_1 & 0 & 0 & \\
A_{-1} & A_0 & A_1 & 0 & \\
0 & A_{-1} & A_0 & A_1 & \ddots  \\
0 & 0 & A_{-1} & A_0 & \ddots  \\
  &   &  \ddots & \ddots & \ddots
\end{bmatrix}.
\]
This is the transition probability matrix among the states in the levels $1$ and above,
avoiding the level $0$.  Since the QBD is assumed to be irreducible, the matrix $W$
defined as $W = \sum_{\nu \geq 0} H^\nu$ converges.   We partition it in blocks $W_{nk}$,
for $n$, $k \geq 1$, and the element $W_{(n,i),(k,j)}$ is the expected number of visits to
 the state $(k,j)$, starting from $(n,i)$, before the first visit to any state in level 0.

\begin{lemma}
   \label{t:w}
If the QBD is irreducible, then
\begin{align}
   \label{e:wnk}
W_{nk} & = G^{n-k}  W_{kk}  \qquad \mbox{for $n \geq k$},  \\
   \label{e:wnkb}
 & = W_{nn}  R^{k-n}  \qquad \mbox{for $n \leq k$},  \\
   \label{e:wkk}
W_{kk} & = \sum_{0 \leq \nu \leq k-1} G^\nu (I-U)^{-1} R^\nu,
\end{align}
for $n$, $k \geq 1$.  Furthermore,
\begin{equation}
   \label{e:sumw}
\sum_{k \geq 1} W_{nk} \vc 1= \vtau_n.
\end{equation}

\end{lemma}

\proof
The matrix $W$ is the minimal nonnegative solution of the system $W (I-H) = I$ and by
Bini {\em et al.}~\cite[Page 102]{blm05},
\[
I-H =
\begin{bmatrix}
I  & -R & 0 \\
0 & I & -R & \ddots \\
 & 0 & I &  \ddots \\
 &  & \ddots & \ddots
\end{bmatrix}
\begin{bmatrix}
I-U & 0 \\
0 & I-U & 0 \\
    & \ddots & \ddots
\end{bmatrix}
\begin{bmatrix}
I & 0 \\
-G & I & 0 \\
 0  & -G & I & \ddots \\
  & \ddots & \ddots & \ddots
\end{bmatrix}
\]
so that
\[
W =
\begin{bmatrix}
I & 0 \\
G & I & 0 \\
G^2 &G & I & 0 \\
\vdots & \vdots & \vdots
\end{bmatrix}
\begin{bmatrix}
(I-U)^{-1} & 0 \\
0 & (I-U)^{-1} & 0 \\
    & \ddots & \ddots
\end{bmatrix}
\begin{bmatrix}
I & R & R^2 & \ldots \\
0 & I & R  & \ldots \\
 & 0 & I   & \ldots \\
 & & & \ddots
\end{bmatrix}
\]
and one proves by direct verification that $W$ given by (\ref{e:wnk}, \ref{e:wnkb}, \ref{e:wkk})
is one solution of $W (I-H)=I$.  To prove that it is minimal, we use the physical meaning of $W$,
and adapt the proof of Latouche {\em et al.}~\cite[Lemma 4.2]{lmt12}.

For $n \geq k$, (\ref{e:wnk}) is proved as follows: starting from level $n$,
the process must first move down to level $k$, with probability matrix $G^{n-k}$, this
justifies the first factor in (\ref{e:wnk}).  The second factor is justified by the fact
that once level $k$ has been reached, the number of visits is given by $W_{kk}$.  For $n \leq k$,
$W_{nn}$ gives the expected number  of visits to level $n$; to each of these visits there corresponds
 a possible excursion to higher levels, and $R^{k-n}$ is the number of visits to level $k$ for each of
 these excursions.  For $n=k$, we count separately the  visits before- and the visits after the first
 visit to level $k-1$; this gives
\begin{align*}
W_{kk} & = (I-U)^{-1} + G W_{k-1,k} \\
 & = (I-U)^{-1} + G W_{k-1,k-1} R \qquad \mbox{by (\ref{e:wnkb}),}
\end{align*}
from which (\ref{e:wkk}) follows.  Finally, (\ref{e:sumw}) immediately results from the definition of $\vtau_n$ and the physical interpretation of $W$.
\BOX

\begin{theorem}\label{t:qbdb}
If the QBD is irreducible, aperiodic and positive recurrent,  then
its deviation matrix is given by
 $\DD = (I - \vone \vpi) K$, where
\begin{align*}
K_{0k} & = (I-P_*)^\# (I - \vtau_0 \vpi_0) R^k & & k \geq 0 \\
K_{n0} & = - \vtau_n \vpi_0 + G^n K_{00} & & n \geq 1 \\
K_{nk} & = W_{nk}  - \vtau_n \vpi_k + G^n K_{0k} & &  n, k \geq 1,
\end{align*}
and $\vtau$ is given in Corollary \ref{t:tau}.
\end{theorem}

\proof
We write the system $(I-P)\DD = I - \vone \vpi$ as follows, so as to make the structure clearly visible:
\begin{equation}
   \label{e:systeme}
\left[
\begin{array}{c|ccc}
I-B &
\begin{array}{ccc}-A_1 & 0 & \ldots \end{array} \\ \hline
-A_{-1} & \\
0 &  I- H\\
\vdots
\end{array}
\right]
\left[
\begin{array}{cccc}
\DD_{00}  & \DD_{01} & \DD_{02} & \ldots \\ \hline
\DD_{10} & \DD_{11} & \DD_{12}  & \ldots \\
\DD_{20} & \DD_{21} & \DD_{22}  & \ldots \\
\vdots & & \vdots
\end{array}
\right]
= \left[
\begin{array}{cccc}
I & 0 & 0 \\ \hline
0 & I \\
0 &  & I \\
 &
\end{array}
\right]
- \vone \vpi.
\end{equation}
We perform one step of Gaussian elimination and isolate level 0 from
the other levels, so that for the levels 1 and above, we obtain
\begin{align}
  \nonumber
&
(I-H)
\begin{bmatrix}
\DD_{10} & \DD_{11} & \DD_{12}  & \ldots \\
\DD_{20} & \DD_{21} & \DD_{22}  & \ldots \\
\vdots & & \vdots
\end{bmatrix}  \\
  \nonumber
& \qquad =
\left[
\begin{array}{cccc}
0 &I & 0 \\
0 & & I \\
\vdots &  & &\ddots
\end{array}
\right]
- \vone \vpi
 +  \vligne{A_{-1} \\ 0 \\ \vdots} \vligne{\DD_{00} & \DD_{01} & \DD_{02} &
   \ldots}
\end{align}
or
\begin{align}
  \nonumber
 \begin{bmatrix}
\DD_{10} & \DD_{11} & \DD_{12}  & \ldots \\
\DD_{20} & \DD_{21} & \DD_{22}  & \ldots \\
\vdots & & \vdots
\end{bmatrix}
  = \  &
W
\left[
\begin{array}{cccc}
0 &I & 0 \\
0 & & I \\
\vdots &  & &\ddots
\end{array}
\right]  - W\vone \vpi
\\
&
 + W \vligne{A_{-1} \\ 0 \\ \vdots} \vligne{\DD_{00} & \DD_{01} & \DD_{02} & \ldots}
   \label{e:systemec}
\end{align}
and for level 0 we obtain
\begin{eqnarray}
  \nonumber
\lefteqn{(I-P_*) \vligne{\DD_{00}  & \DD_{01} & \DD_{02} & \ldots} =
\vligne{I & 0 & 0 & \ldots} - \vone \vpi}  \qquad \qquad \\
 & + &\vligne{A_1 & 0 & 0 & \ldots} W
\big(
\left[
\begin{array}{cccc}
0 &I &  \\
0 &  & I \\
\vdots &  & & \ddots
\end{array}
\right]
- \vone \vpi \big)
   \label{e:doo}
\end{eqnarray}
as the matrix $N_B$ in (\ref{e:wam1}) is the matrix $W$ here.
By (\ref{e:wkk}), and the fact that  $R=A_1 (I-U)^{-1}$, we have
$\vligne{A_1 & 0 & 0 & \ldots} W = \vligne{ R & R^2 & R^3 & \ldots}$,
and (\ref{e:doo}) becomes
\begin{align*}
(I-P_*) & \vligne{\DD_{00}  & \DD_{01} & \DD_{02} & \ldots} \\
& = \vligne{ I & R & R^2  & \ldots} - \sum_{\nu \geq 0} R^\nu \vone \vpi \\
& = \vligne{ I & R & R^2  & \ldots} - (I-R)^{-1} \vone \vpi_0\vligne{ I & R & R^2  & \ldots}  \\
& = (I - \vtau_0 \vpi_0) \vligne{ I & R & R^2  & \ldots}
\end{align*}
by (\ref{3.3}).  By (\ref{e:drazin}), this gives us
\begin{equation}
   \label{e:dok}
\DD_{0k} = (I-P_*)^\# (I-\vtau_0 \vpi_0) R^k + \vone \valpha_k = K_{0k} + \vone \valpha_k,
\qquad \mbox{for $k \geq 0$,}
\end{equation}
where the vectors $\valpha_k$ will be determined later.

Returning to (\ref{e:systemec}), we find that
\begin{align}
  \nonumber
& \begin{bmatrix}
\DD_{10} & \DD_{11} & \DD_{12}  & \ldots \\
\DD_{20} & \DD_{21} & \DD_{22}  & \ldots \\
\vdots & & \vdots
\end{bmatrix}  \\
   \label{e:systemeb}
& \qquad  = \vligne{\vzero & W} - \vligne{\vtau_1 \\ \vtau_2 \\ \vdots} \vpi + \vligne{G \\ G^2 \\
\vdots}\vligne{\DD_{00} & \DD_{01} & \ldots}
\end{align}
by (\ref{e:wam1}, \ref{e:sumw}).  Since the QBD is positive recurrent, $G$ is stochastic and,
combining (\ref{e:dok}, \ref{e:systemeb}), we find that $\DD_{nk} = K_{nk} + \vone \valpha_k$
for all $n$ and $k$, or $\DD = K + \vone \valpha$ in global form.  Since $\DD$ must satisfy
the constraint $\vpi \DD = \vzero$, this shows that $\valpha = -\vpi K$, which concludes the proof.
\BOX

\section{Perturbation analysis}
   \label{s:perturbation}

We  investigate in this section the link between Poisson's equation
and perturbation analysis.  We do not restrict ourselves to QBDs.
Instead, we consider a general Markov chain, which is assumed to be
$w$-geometrically ergodic.
 Recall that an irreducible, aperiodic and positive recurrent QBD is $w$-geometrically ergodic,
 with drift condition given in Lemma~\ref{t:drift}.

Let $P(\delta) = P + \delta Q $ be an irreducible and stochastic
transition matrix, where $Q \vc 1 = 0$ and $\delta$ belongs to a
neighborhood of 0. We may interpret $P(\delta)$ as a perturbation of
$P$. Suppose that $P(\delta)$ is also positive recurrent with
invariant distributions $\boldsymbol{\pi}(\delta)$, for $\delta$
small enough,  and define
$\omega(\delta)= \boldsymbol{\pi} (\delta) \boldsymbol{g}$.
We are interested in  the derivatives of $\vpi(\delta)$ and
$\omega(\delta)$ with respect to $\delta$, evaluated at $\delta =0$.
Recall that the drift condition
\textbf{D}($\boldsymbol{v}, \lambda, b, C$) was introduced at the end of
Section~\ref{s:basicPoisson},

\begin{proposition}
   \label{t:perturb}
Let $\{X_t\}$ be an irreducible, aperiodic and positive recurrent Markov
chain. Assume that $P$ satisfies the drift condition
\textbf{D}($\boldsymbol{v}, \lambda, b, C$), that
$\|Q\|_{\boldsymbol{v}}<\infty$, and that $\| \vc g \|_{v} <
\infty$. One has
\begin{equation}\label{deri-pi}
   \left.\frac{\ud^n \boldsymbol{\pi}}{\ud \delta^n}\right|_{\delta =0}= n! \boldsymbol{\pi} (Q \mathcal{D})^n,
\end{equation}
and
\begin{equation}\label{deri-delta}
    \left.\frac{\ud^n \omega}{\ud \delta^n}\right|_{\delta =0} = n! \boldsymbol{\pi} (Q \mathcal{D})^n \boldsymbol{g}.
\end{equation}
In particular, for $n=1$,
\begin{equation}\label{deri-po}
  \left.\frac{\ud \omega}{\ud \delta}\right|_{\delta =0} =
  \boldsymbol{\pi} Q \boldsymbol{h},
\end{equation}
 where $\boldsymbol{h} = \DD \vc g$ is a solution of Poisson's equation such that $\boldsymbol{\pi} \vc h = 0$.
  \end{proposition}

\proof Since  $P$ satisfies \textbf{D}($\boldsymbol{v}, \lambda, b, C$), it is ${v}$-geometrically ergodic and so
\[
  \|\mathcal{D}\|_{{v}} \leq \sum_{n=0}^\infty \|P^n-\vone\vpi\|_{{v}} \leq \sum_{n=0}^\infty  r \beta^{n}<\infty.
\]
Define $q = \|Q\|_v$ and take $\delta<{1-\lambda}/{q}$.  We have
\[
P(\delta) \boldsymbol{v} = (P+\delta Q) \boldsymbol{v} \leq (\lambda + \delta q)\boldsymbol{v} + b \indic[C]
\]
where $\lambda + \delta q<1$. Thus, $P(\delta)$ is also ${v}$-geometrically ergodic and
its invariant probability measure $\boldsymbol{\pi}(\delta)$ exists.
We pre-multiply both sides of $(I-P)\DD = I - \vone \vpi$ by $\vpi(\delta)$,
and use $\boldsymbol{\pi}(\delta) P(\delta) = \boldsymbol{\pi}(\delta)$, to obtain
\[
  \boldsymbol{\pi}(\delta)-\boldsymbol{\pi} = \boldsymbol{\pi}(\delta) \mathcal{D}- \boldsymbol{\pi}(\delta)
  P \mathcal{D}= \boldsymbol{\pi}(\delta) P(\delta) \mathcal{D}- \boldsymbol{\pi}(\delta)
  P \mathcal{D}=\delta \boldsymbol{\pi({\delta})} Q \mathcal{D}.
\]
If $\delta$ is small enough, so that $\delta<(1-\lambda)/q$ and $\delta \| Q \mathcal{D}\|_{\boldsymbol{v}} <1$, then
\[
   \boldsymbol{\pi}(\delta)=\boldsymbol{\pi} (I- \delta Q \mathcal{D})^{-1} =\boldsymbol{\pi}  \sum_{n=0}^{\infty}
    (Q \mathcal{D})^n \delta^n,
\]
and, since the power series is convergent, we obtain  (\ref{deri-pi}).
 If  $\|\vc g \|_v < \infty$, then
$
   \omega(\delta) =\sum_{n\geq 0} \boldsymbol{\pi} (Q \mathcal{D})^n \vc g \delta^n,
$
from which (\ref{deri-delta}) directly  follows.   In particular,
$\ud \omega/\ud \delta |_{\delta = 0}= \vpi Q \DD \vc g$.
 \BOX

\begin{remark}  \rm
This property is proved in \cite{cc97}  under the assumption that
the Markov chain is uniformly (or strongly) ergodic. Generally
speaking, strong ergodicity is  stricter than one might wish.
Indeed, we see from  Proposition 2.1 in \cite{hl04} that many
discrete-time Markov chains are not strongly ergodic, since their
transition matrix is a Feller transition matrix, that is,
$\lim_{i\rightarrow \infty} P_{ij}=0$ for any fixed $j$. Proposition
\ref{t:perturb} here is an improvement since it requires  the weaker
condition of geometric ergodicity.

We also note that the same assumption of uniform ergodicity is made in
Altman {\em et al.}~\cite{aan04} and in Liu~\cite{l12}, with the added
constraint that the set $C$ is a single state. In our drift condition,
$C$ may be a finite set which, as we have seen, is more convenient
for matrix-analytic models.
Performance analysis of Markov chains on a general state space are
analyzed in Kartashov~\cite{karta86} and Heidergott and
Hordijk~\cite{hh03}, under conditions which are essentially equivalent
to geometric ergodicity.

\end{remark}

\section{Application to a queue}
   \label{s:exemple}

To illustrate our results, we consider the PH/M/1 queue.  This is a
system with a single server and a buffer of unlimited capacity, the
service times distribution is exponential, with parameter $\mu$, the
arrivals form a renewal process, and the intervals of time between
arrivals have a PH distribution, with representation $(\vc\sigma,
S)$.  These queues are {\em continuous-time} QBDs with generator
\[
Q =
\begin{bmatrix}
S & \vc s \vc\sigma & \\
\mu I & S - \mu I & \vc s \vc\sigma & \\
 & \mu I & S - \mu I & \ddots \\
 & & \ddots & \ddots
\end{bmatrix}.
\]
Their stationary expected queue length $L$ is easily seen to be equal to
\[
L = \sum_{n \geq 0} n \vpi_0 R^n \vone = \vpi_0 R (I-R)^{-2} \vone.
\]

As our results have been formulated for discrete-time Markov chains,
we uniformize the PH/M/1 queue and obtain the transition
matrix~(\ref{e:P}) with $A_{-1}= (\mu/\gamma) I$, $A_0 =
I+(1/\gamma)(S-\mu I)$, $A_1 = (1/\gamma) \vc s \vc\sigma$, and $B =
A_0+A_{-1}$, where $\gamma > \mu + S_{ii}$ for all $i$.
We denote the process  as $\{(Y_n, \varphi_n)\}$, where $Y_n$ is the level at
time $n$ and $\varphi_n$ is the phase, and we recall that the
uniformized QBD has the same stationary distribution as the PH/M/1
queue itself.

We are interested in the vector $\vc m$ with
\begin{equation}
   \label{e:m}
m_{\ell,j} = \sum_{n \geq 0} (\E[Y_n|Y_0=\ell,\varphi_0=j] / L -1).
\end{equation}
Coolen-Schrijner and van Doorn~\cite{cv02} discuss the
similar quantity
\[
m({\cal Y}) = \int_0^\infty (1-\E[Y(t)]/\E[Y]) \, \ud t,
\]
where $\{Y(t)\}$ is a continuous-time, stochastically monotone
Markov chain, and $\E[Y]$ is its stationary expectation.  For such a
Markov chain starting from the minimal state, $\E[Y(t)] < \E[Y]$ and
the finite-time expectation monotonically converges to its limit, so
that $m({\cal Y})$ is positive, and it may be interpreted as a
measure of the speed of convergence to stationarity of $\E[Y(t)]$.
This quantity is expressed in \cite{cv02} in terms of the deviation
matrix.

QBD processes are usually not stochastically monotone.
Nevertheless,  a quantity such as $\vc m$  defined in
(\ref{e:m}) is
interesting because it   measures the sensitivity of the expected queue size to the
initial state.
As we show in the following proposition, this monotonicity condition
may be replaced by geometric ergodicity  ---  we also adapt the
formulation of the property to the case of discrete-time Markov
chains.

\begin{proposition}
   \label{t:speed}
   Let $\{X_t\}$ be an irreducible, aperiodic and positive recurrent
   discrete-time Markov chain, assume that its transition matrix
   satisfies the drift condition \textbf{D}($\boldsymbol{v}, \lambda,
   b, C$). Define
\[
m_i(g) = \sum_{t \geq 0} (\E[g_{X_t}|X_0=i] / \vpi \vc g -1),
\]
where $\vc g$ is a vector such that $\vpi |\vc g| < \infty$ and $\|\vc
g\|_v < \infty$.

One has $\vc m = (\vpi \vc g)^{-1} \DD \vc g$, where $\DD$ is the deviation matrix of the Markov chain.
\end{proposition}

\proof Since $\{X_t\}$ satisfies \textbf{D}($\boldsymbol{v}, \lambda, b,
C$), it is ${v}$-geometrically ergodic,
 and there exist positive constants $r$ and $\beta<1$ such that
\[
  | \E[g_{X_t}|X_0=i]- {\boldsymbol \pi \vc g}|\leq r v_i \|\vc g\|_v \beta^t
\]
for all $i$. Thus,
\[
   \sum_{t \geq 0} \left| \E[g_{X_t}|X_0=i]-E[g_{X_{\infty}}]\right|=\sum_{t \geq 0} \sum_{j\in \mathbb{E}}
   |(P^t_{ij}-{\pi_j}) g_j|<\infty,
\]
which implies
\[
  \sum_{t \geq 0} \left(\E[g_{X_t}|X_0=i]-E[g_{X_{\infty}}]\right) = \sum_{j\in\mathbb{ E}} \sum_{t \geq 0}
  (P^t_{ij}-{\pi_j}) g_j = \sum_{j\in \mathbb{E}} {\cal D}_{ij} g_j=({\cal D} g)_i.
\]
If $\vpi \vc g$ is finite, we divide both sides by $\vpi \vc g$, and the proof is complete.
\BOX

We take $\vc g_n = n \vone$ for all $n$, so that the vector $\vc m$
defined in (\ref{e:m}) is identical to $\vc m(g)$, and measures,
therefore, the difference over the history of the process between the
time-dependent expected queue length and its stationary value $L =
\vpi \vc g$.

By Proposition~\ref{t:speed}, $\vc m = L^{-1} \DD \vc g$, so that $\vc m$ is the solution of
\[
(I-P) \vc m = L^{-1} (I - \vone \, \vpi) \vc g =  L^{-1} \vc g - 1,
\]
with the added constraint that $\vpi \vc m=0$.  We apply Theorem~\ref{t:qbda} and obtain
\begin{align}
   \label{e:mo}
\vc m_0 & = -(\vc\sigma \vc y_1)(S + \vc s \, \vc\sigma G)^\# \vc s   + c_0 \vone \\
   \label{e:mn}
\vc m_n & = \vc y_n + G^n \vc m_0
\end{align}
for all $n \geq 1$, where
\begin{align}
  \nonumber
\vc y_n & = L^{-1} \sum_{0 \leq i \leq n-1}  (n-i) G^i \vtau_1 - \sum_{0 \leq i \leq n-1} G^i\vtau_1 \\
   \label{e:yn}
& \quad  +  L^{-1} \sum_{0 \leq i \leq n-1} G^i  (I-U)^{-1}(I-R)^{-1} A_1\vtau_1
  \\
   \label{e:tauone}
\vtau_1 & = (I-U)^{-1}(I-R)^{-1} \vone
\end{align}
and $c_0$ is such that $\vpi \vc m=0$.  The detailed proof, and the
expression for $c_0$, are not enlightening and they are given in
appendix.

We have considered three different distributions for the inter-arrival
times:
\begin{itemize}
\item Erlang with
 \[
\vc\sigma = \vligne{1 & 0}, \qquad
S = \left[\begin{array}{rr}  -2 & 2 \\ 0 & -2 \end{array} \right];
\]
\item Exponential with $\vc \sigma = \vligne{1}$, $S =
  \vligne{-1}$;
\item Hyper-exponential with
\[
\vc\sigma = \vligne{0.11270167 & 0.88729833}, \qquad
S = \vligne{-0.225403332 & 0 \\ 0 & -1.77459667}
\]
(from Ramaswami and Latouche \cite[page 646]{rl89b}).
\end{itemize}
In all cases, the intervals between successive arrivals have
expectation equal to 1 and the traffic coefficient $\rho$ is equal to
$1/\mu$.  The variance of the Erlang distribution is equal to 0.5 and
that of the Hyper-exponential distribution is~4.

\begin{figure}
\centering
\includegraphics[scale=0.61]{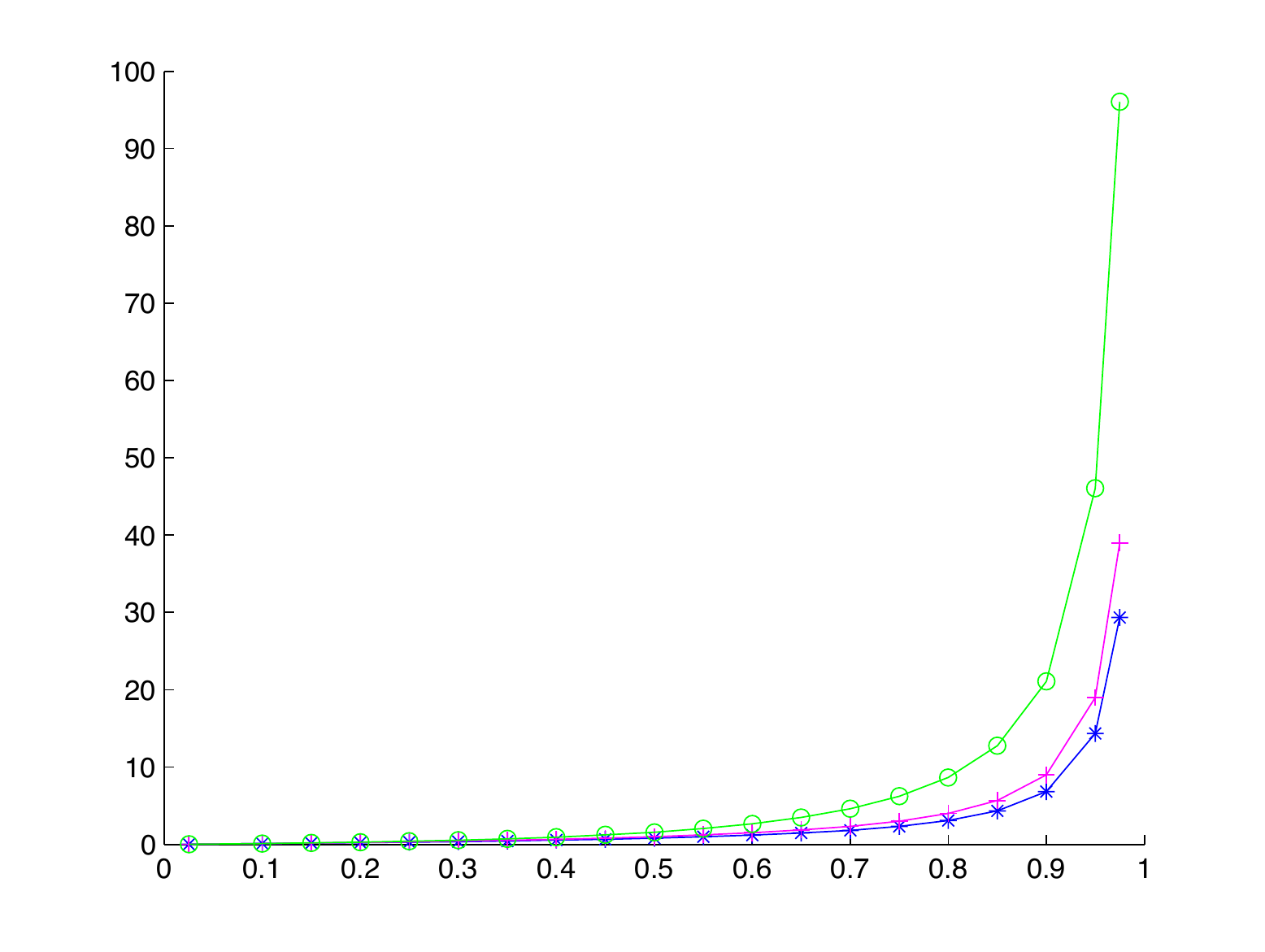}
\caption{Expected queue length as a function of $\rho$.  The curve for
  the E$_2$/M/1 queue is marked with $*$, for the M/M/1 queue with $+$ and for the H$_2$/M/1 with $\circ$. }
   \label{f:mean}
\end{figure}

\begin{figure}
\vspace{-5\baselineskip}
\centering
\includegraphics[scale=0.61]{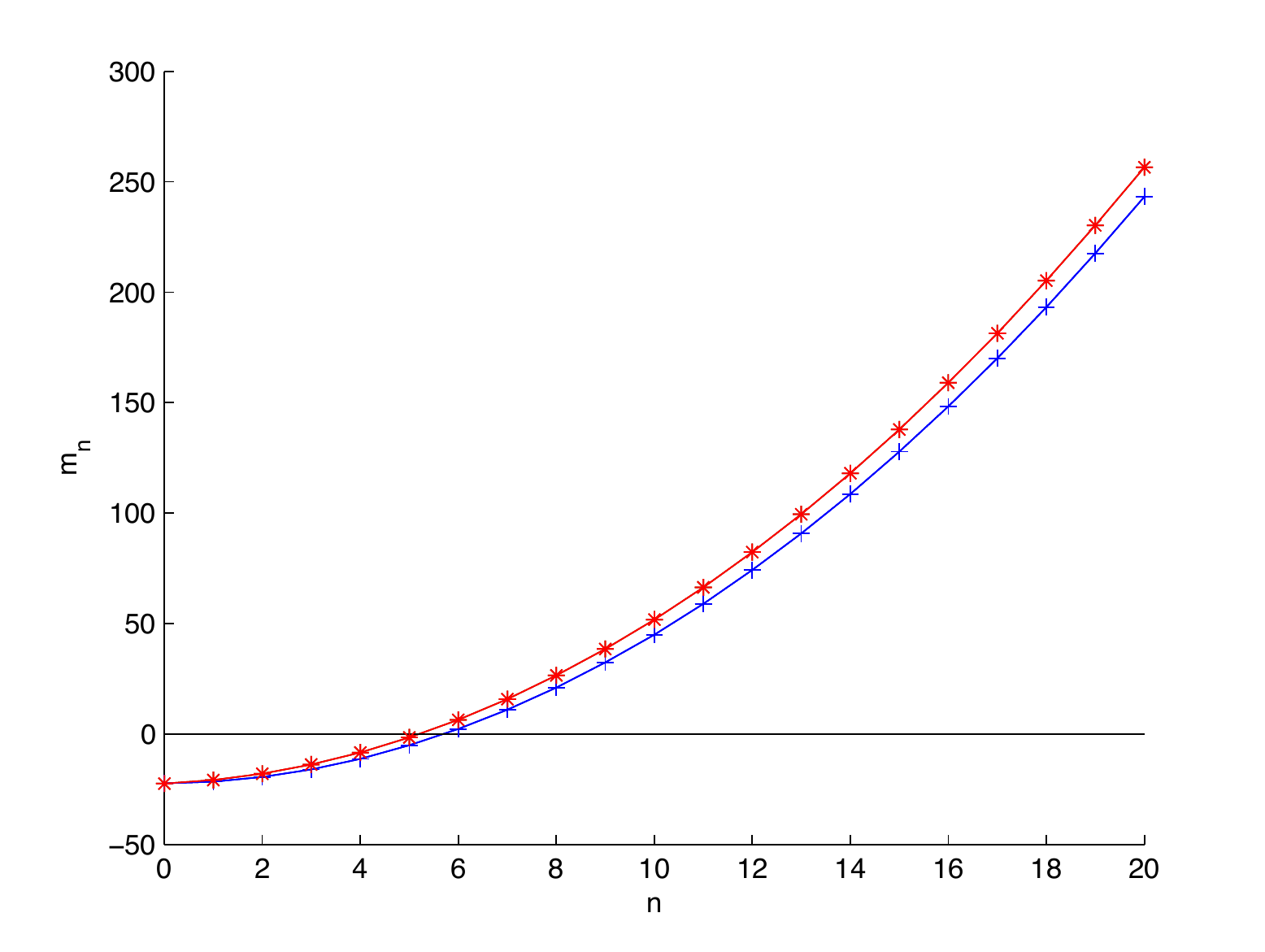}

\includegraphics[scale=0.61]{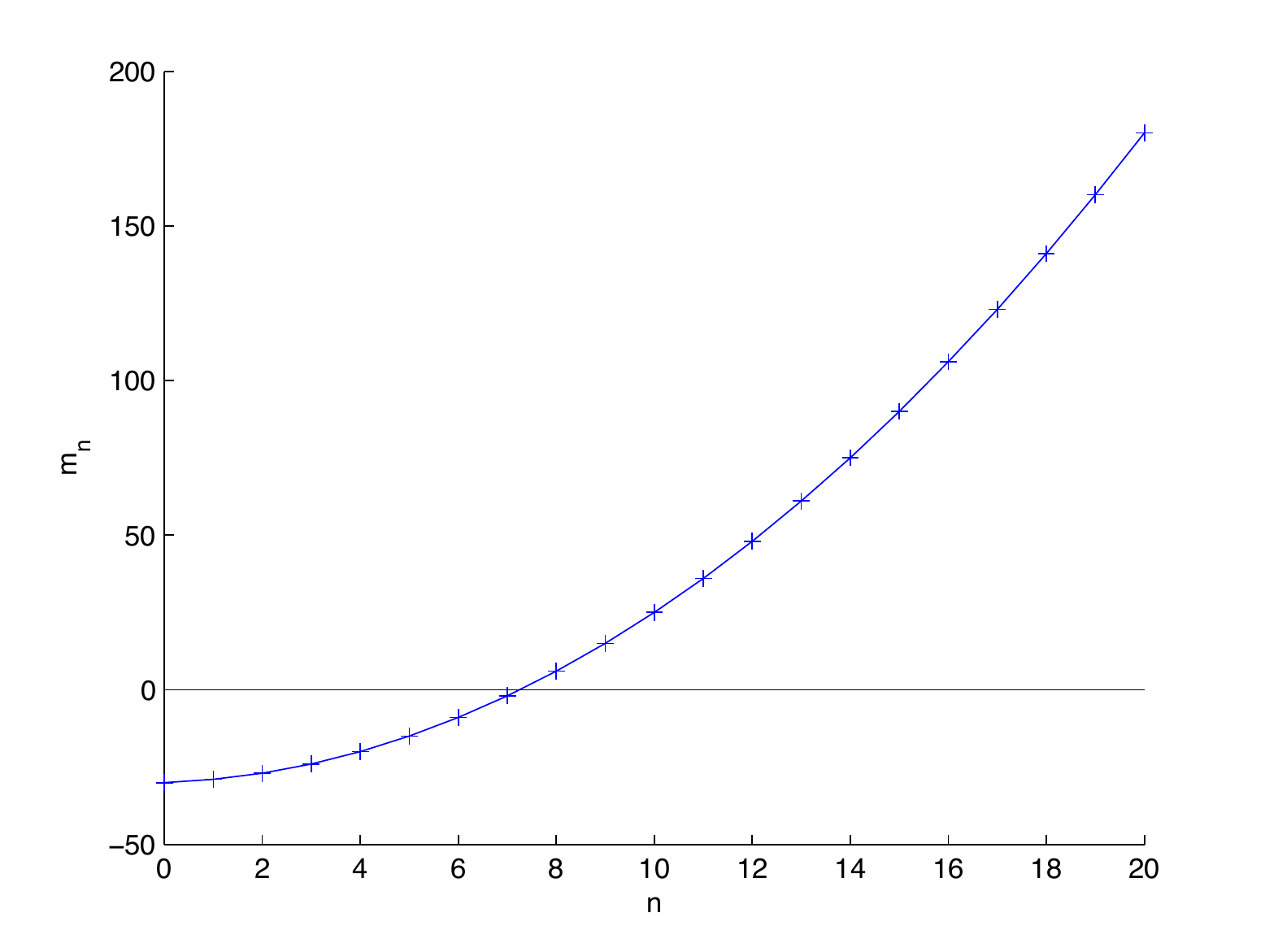}

\includegraphics[scale=0.61]{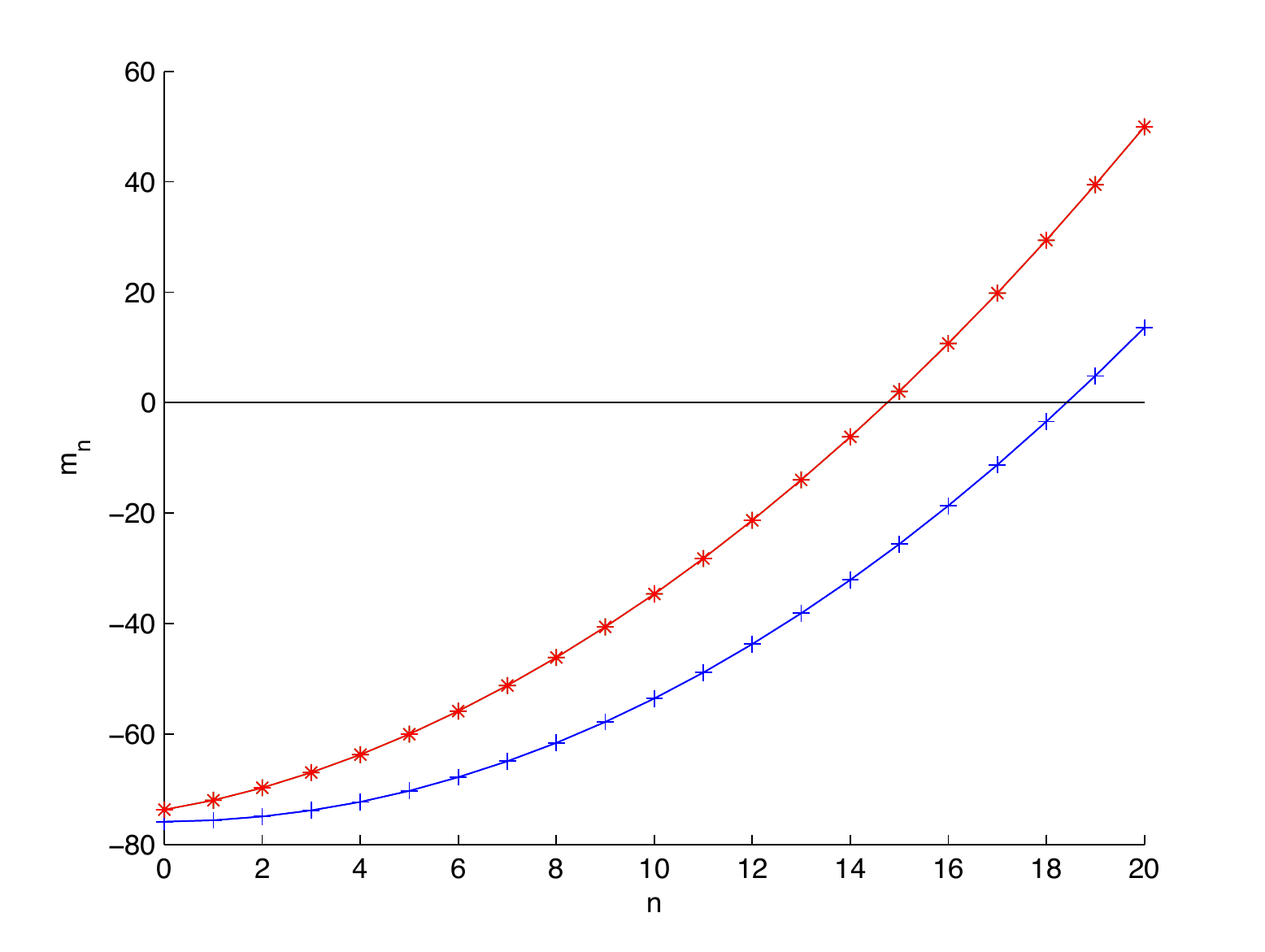}
\caption{Queues E$_2$/M/1 (top), M/M/1 (middle) and H$_2$/M/1
  (bottom);  display of $ m_{n,i}$ as a function of $n$ for each phase.  The traffic
  coefficient is  $\rho= 1/1.2$. The blue curve (marked with ``+'') is for phase 1, the red curve
  (marked with ``*'') is for phase 2. } \label{f:speedb}
\end{figure}

The difference in variability is reflected in the expected stationary
queue length, as shown in Figure~\ref{f:mean} where we plot the value
of $L$ as a function of $\rho$ for the three queues.   We see on
Figure~\ref{f:speedb} that it is also reflected in the sensitivity.

We take $\mu = 1.2$ fixed, so that $L$ takes respectively the values
3.8, 5, and 11.1 for the Erlang, the exponential and the
hyper-exponential arrival processes.  We show on Figure~\ref{f:speedb}
the components of the vectors $\vc m_n$ for small values of $n$.  In
all cases, the difference is increasing, starting from a negative
value for small values of $n$, and becoming positive at some $n$
greater than $L$.  We also observe a clear difference between the
plots for the E$_2$/M/1 and the H$_2$/M/1 queues.  Not surprisingly,
the influence of the first phase is more pronounced in the latter
case.  Finally, the range of values for $\vc m_n$, $0 \leq n \leq 20,$
is greater when the arrival distribution has a smaller variance, a
feature that we did not expect.

\subsection*{Acknowledgment}

It is a pleasant duty to mention that Bernd Heidergott  and three
anonymous referees have offered numerous and insightful comments on
the first draft of this paper.

The authors thank the Minist\`ere de la Communaut\'e fran\c{c}aise de
Belgique for funding this research through the ARC grant
AUWB-08/13--ULB~5.

The third author also gratefully acknowledges the
support of the Fundamental Research Funds for the Central Universities
(grant number 2010QYZD001) and the National Natural Science Foundation
of China (grant number 10901164).

\appendix
\section*{Appendix}

For the PH/M/1 queue, (\ref{3.5}) becomes
\begin{align*}
\vc y_0 &= L^{-1} \vc g_0 - \vone + A_1 \vc y_1 \\
  & = -\vone + (1/\gamma)  (\vc\sigma \vc y_1)  \vc s
\end{align*}%
and (\ref{3.1}) gives
\begin{align*}
\vc m_0 & = (I - P_*)^\# \vc y_0 + c_0 \vone \\
  & = - \gamma (S + \vc s \, \vc\sigma G)^\# \vc y_0 + c_0 \vone \\
  & = - (\vc\sigma \vc y_1) (S + \vc s \, \vc\sigma G)^\#  \vc s  + c_0 \vone
\end{align*}
since $(S + \vc s \, \vc\sigma G)^\# \vone = \vzero$, which proves (\ref{e:mo}).

Equation (\ref{e:mn}) is merely a reformulation of (\ref{3.6}), so that we only need to prove (\ref{e:yn}).
 Before doing so, we observe that
\[
(I-U)^{-1} \sum_{n \geq 0} R^n \vone = (I-U)^{-1}(I-R)^{-1} \vone
\]
and is equal to $\vtau_1$ by (\ref{3.4}) --- this proves (\ref{e:tauone}) --- and we also note that
\[
\sum_{n \geq 0} n R^n \vone = (I-R)^{-1} R (I-R)^{-1} \vone = (I-R)^{-1} A_1 \vtau_1
 \]
since $R = A_1 (I-U)^{-1}$.
By (\ref{3.6}),
\begin{align*}
\vc y_n & = \sum_{0 \leq i \leq n-1} G^i (I-U)^{-1} \sum_{k \geq 0} R^k (L^{-1} (n-i+k) \vone - \vone) \\
 & = L^{-1} \sum_{0 \leq i \leq n-1}  (n-i) G^i \vtau_1 -   \sum_{0 \leq i \leq  n-1} G^i \vtau_1  \\
  & \quad +      L^{-1}\sum_{0 \leq i \leq n-1} G^i (I-U)^{-1} \sum_{k \geq 0} k R^k \vone
\end{align*}
from which (\ref{e:yn}) follows.

Finally, we need to determine $c_0$ such that $\vpi \vc m = 0$.  Observe that $\sum_{n \geq 0} \vpi_n G^n = \vpi_0 M$,
 where
$
M = \sum_{n \geq 0} R^n G^n
$
\ is the unique solution of the linear system $M = I + R M G$ and may easily be computed.
We immediately obtain from (\ref{e:mo}, \ref{e:mn}) that
\begin{equation}
   \label{e:c}
c_0  =  (\vc\sigma \vc y_1) \vpi_0 M (S + \vc s \, \vc\sigma G)^\# \vc s   - \sum_{n \geq 1} \vpi_0 R^n \vc y_n.
\end{equation}
To evaluate the last term, we proceed in three steps.  Firstly, we write
\begin{align*}
\sum_{n \geq 1} R^n \sum_{0 \leq i \leq n-1} G^i & = \sum_{i \geq 0}  \sum_{n \geq i+1} R^n G^i \\
 & = \sum_{i \geq 0}  R^i \sum_{n \geq 1} R^n G^i \\
 & = \sum_{i \geq 0}  R^i R (I-R)^{-1} G^i \\
 & = R (I-R)^{-1} M.
\end{align*}
Secondly,
\begin{align*}
\sum_{n \geq 1} R^n \sum_{0 \leq i \leq n-1} (n-i) G^i & = \sum_{i \geq 0} R^i \sum_{n \geq i+1} (n-i) R^{n-i} G^i \\
 & = \sum_{i \geq 0} R^i \sum_{n \geq 1} n R^{n} G^i \\
 & = R (I-R)^{-2} M.
\end{align*}
Finally, we use (\ref{e:yn}, \ref{e:c}) and find that
\begin{align}
   \nonumber
c_0  & = (\vc\sigma \vc y_1) \vpi_0 M (S + \vc s \, \vc\sigma G)^\# \vc s
  - L^{-1} \vpi_0 R (I-R)^{-2} M \vtau_1 \\
  \nonumber
& \quad + \vpi_0 R (I-R)^{-1} M \vtau_1 \\
   \label{e:cc}
 & \quad - L^{-1} \vpi_0 R (I-R)^{-1} M (I-U)^{-1} (I-R)^{-1} A_1 \vtau_1.
\end{align}


\end{document}